\input amstex
\mag=\magstep1
\documentstyle{amsppt}
\NoBlackBoxes
\def\fbox#1#2{\vbox to 0pt{\hrule\hbox{\hsize=#1\vrule\kern2pt
  \vbox{\kern2pt#2\kern2pt}\kern2pt\vrule}\hrule}}

\def \-#1{{#1}^{-1}}

\def \t#1{\tilde{#1}}
\def \+n#1{{#1}^{+\nu}}
\def \o#1{\overline{#1}}
\def \ex#1#2#3#4#5#6
{0\to \Cal O_{#1}(#2) \to \Cal O_{#3}(#4) \to \Cal O_{#5}(#6) \to 0}


\topmatter
\title 
The existence of Gorenstein Terminal Fourfold Flips
\endtitle
\author Hiromichi Takagi$^*$
\endauthor
\address Department of mathematical science, 
Tokyo University ,Komaba,Tokyo,Japan 153
\endaddress
\thanks 1991 Mathematics Subject Classification. Primary: 14E30, 
Secondary: 14J35, 14M05.\newline 
\ \ $*$Research Fellow of the Japan Society for the Promotion of Science.
\endthanks
\abstract
We prove that 
for a flipping contraction from a Gorenstein terminal $4$-fold, 
the pull back of a general hyperplane section of the down-stair 
has only canonical singularities. 
Based on this fact and using Siu-Kawamata-Nakayama's
extension theorem [Si], [Kaw4], [Kaw5] and [Nak2], we prove the existence 
of the flip of such a flipping contraction.
Furthermore we classify such flipping contractions and the flips
under some additional assumptions.
\endabstract
\email htakagi\@ ms.u-tokyo.ac.jp 
\endemail
\endtopmatter     

\head 0.Introduction
\endhead
To proceed the Minimal Model Program (in short MMP), 
an elementary transformation called
flip is very important (see [KMM] for detail). 
\definition{Definition 0.1}
Let $X$ be a normal algebraic variety (resp. normal analytic variety)
with only canonical singularities
and $Y$ a normal algebraic variety 
(resp. $(Y,S)$ a pair of an analytic space and its compact subspace). 
A projective morphism $f:X\to Y$ is called a flipping contraction if
\roster
\item $-K_X$ is $f$-ample;
\item $\rho(X/Y) = 1$ (resp. $\rho(X/Y, f^{-1}(S))=1$);
\item $f$ is an isomorphism in codimension $1$.
\endroster
If there exists a normal algebraic variety (resp. normal analytic variety) 
$X^+$ 
with only canonical singularities and a projective morphism
$f^+:X^+\to Y$ such that
\roster
\item $K_{X^+}$ is $f^+$-ample;
\item $f^+$ is an isomorphism in codimension $1$,
\endroster
we call $f^+$ the flip of $f$.
We call the following diagram a flipping diagram: 
$$\matrix
X & \dashrightarrow & {X^+} \\
{f\searrow} & \ &  {\swarrow f^+} \\
 \ & Y & \ & .
\endmatrix $$
\enddefinition

The existence of the flip is a very hard problem.
In dimension $3$, Shigefumi Mori proved it in [M4]. 
As a test case of $4$-dimensional flips, we consider a flipping
contraction from an algebraic $4$-fold with only Gorenstein terminal
singularities.
Let $X$ be an algebraic $4$-fold with only Gorenstein terminal
singularities and $f\: X\to Y$ be a flipping contraction.
Let $E$ be the exceptional locus.
Since there is no flipping contraction from an algebraic (or analytic)
$3$-fold with only Gorenstein terminal singularities,
we find that $f(E)$ is a set of finite points.
Hence replacing $Y$ by a small Stein neighborhood of a point in $f(E)$,
we can proceed in the analytic category. Precisely speaking,
we consider the following object below (we call this $(*)$).

$(*)$
Let $X$ be an analytic $4$-fold with only Gorenstein terminal
singularities and $(Y,P)$ a pair of a contractible $4$-dimensional Stein space
and a point in it
such that $Y$ has only ccDV singularities (i.e., singularities whose general
hyperplane sections have only cDV singularities) outside $P$.
Let $f\: X\to Y$ be a flipping contraction and $E:=f^{-1}(P)$, 
i.e., the exceptional locus of $f$.

In [Kaw3], 
Yujiro Kawamata considered the case where $X$ is
smooth. He proved the following:

\proclaim{Theorem 0.1}
Assume that $X$ is smooth.
Then  the flip exists and 
$E\simeq \Bbb P^2$ and
$\Cal N_{E/X}\simeq \Cal O_{\Bbb P^2}(-1)\oplus \Cal O_{\Bbb P^2}(-1)$. 
In particular we obtain the flip by blowing up $E$ 
(the exceptional locus of the blowing up is $\Bbb P^2\times \Bbb P^1$)
and
blowing down this $\Bbb P^2\times \Bbb P^1$ to $\Bbb P^1$. 
\endproclaim

Quite recently Yasuyuki Kachi proved in his preprint [Kac2] 
the following:

\proclaim{Theorem 0.2}
Assume that $X$ is singular and
has only isolated complete intersection terminal singularities. 
Suppose that there is a member of $|-2K_Y|$ through $P$ 
which has only a rational singularity at $P$. 

Then the flip exists and $E\simeq \Bbb P^2$ and
$\Cal N_{E/X}\simeq \Cal O_{\Bbb P^2}\oplus \Cal O_{\Bbb P^2}(-2)$.
Furthermore $X$ has only one singularity on $E$, which is analytically
isomorphic to $o\in (xy + zw + t^m = 0) \subset \Bbb C^4$.  
\endproclaim
He proved the existence of such a flip by induction
and constructed the desired flip very explicitly 
(see [Kac2, $\S 8$] for detail).
He also investigated some special semistable $4$-fold flipping contractions
in [Kac1].

Our starting points are the following two theorems: 

\proclaim {Theorem 1.2 (Rough classification of the exceptional locus)} 
Assume that  
the exceptional locus $E$ contains $2$-dimensional components. 
Let $E=\cup E_i$ be  the irreducible decomposition of $E$. 
Then 
$E$ is purely $2$-dimensional and $(E_i, -K_X|_{E_i})$ is isomorphic to 
$(\Bbb F_{n,0}, nl)$, where $l$ is a ruling of $\Bbb F_{n,0}$.
\endproclaim

\proclaim{Theorem 1.3}
Let $B$ a general hyperplane section through $P$.
Then the strict transform $A:=f^*B$ has only canonical singularities.
\endproclaim
Our main result is the following:

\proclaim{Main Theorem (See Corollary 2.2)} 
The flip of $f$ as in $(*)$ exists.
\endproclaim
We will prove the finite generation of
$\oplus _{m\geq 0} f_* \Cal O (mlK_X)$ ($l$ is the index of $X$)
using Theorem 1.3 and
the Siu-Kawamata-Nakayama's extension theorem (see Theorem 1.5).

By using the existence of the flip, we obtain a rough classification
of $f$ and $f^+$ (see Corollary 2.3) and 
Furthermore if $A$ (as in Theorem 1.3) 
has only isolated singularities and $E$ is irreducible,
We can give the more detailed description (see Corollary 2.4).

\definition{Acknowledgment}
This is a revised version of my paper [T].
I express my hearty thanks to Professor Yujiro Kawamata for giving me 
useful comments and encouraging me.
I am grateful to Professor Keiji Oguiso for listening to my talk on this
subject carefully. 
I am thankful to Doctor Yasuyuki Kachi for stimulus discussions on this 
subject.
I also thank Professor Shigeru Mukai and Professor Eiichi Sato
for giving me the opportunity to talk about this paper at Kyushu University.
Finally I thank the referee for pointing out my English mistakes and
insufficient exposition in the proof of Claim 2.6 (1).
\enddefinition

\definition{Notation and Convention}
\roster
\item In this paper, we will work over $\Bbb C$, the complex number
field and in the analytic category;
\item We denote by $\Bbb F_n$, the Hirzebruch surface 
$\Bbb P(\Cal O_{\Bbb P^1} \oplus \Cal O_{\Bbb P^1}(-n))$ 
and by $\Bbb F_{n,0}$ the normal surface 
which is obtained from the Hirzebruch surface $\Bbb F_n$ by contracting 
the negative section. 
\endroster
\enddefinition

\head 1. Preliminaries \endhead

\proclaim{Theorem 1.1}
Let $X$ and $Y$ be normal log terminal varieties and $f: X\to Y$
a projective morphism.
Let $L$ a $f$-ample line bundle on $X$ and $F$ a fiber of $f$. 
Assume that $f\:X\to Y$ is the adjoint 
contraction supported by $K_X+rL$
and either $\dim F<r+1$ if $\dim Y<\dim X$ or $\dim F\leq r+1$ if
$\dim Y= \dim X$.

Then $f^*f_*L\to L$ is surjective at every point of $F$.
\endproclaim
\demo{Proof}
 See [AW1]. They assume that $L$ is ample but their proof works also for the
case that $X$ is analytic and $L$ is relatively ample.
\qed
\enddemo

\proclaim {Theorem 1.2 (Rough classification of the exceptional locus)} 
We consider the object $(*)$.
Assume that 
the exceptional locus $E$ contains $2$-dimensional components. 
Let $E=\cup E_i$ be  the irreducible decomposition of $E$. 
Then 
$E$ is purely $2$-dimensional and $(E_i, -K_X|_{E_i})$ is isomorphic to 
$(\Bbb F_{n,0}, nl)$, where $l$ is a ruling of $\Bbb F_{n,0}$.
\comment
\item $E$ has at most $3$ components.
If $E$ has more than $2$ components, $E_i$ intersects one smooth rational curve
which is a ruling on each $E_i$.      
\endcomment
\endproclaim

\demo{Proof}
By Theorem 1.10 and 1.19 of [AW2], it is sufficient to exclude the 
following possibilities:
$(E_i, -K_X|_{E_i})$ is isomorphic to $(\Bbb P^2, \Cal O_{\Bbb P^2}(2))$
or $(\Bbb F_n, C_0+ml)$,
where $C_0$ is the negative section and $l$ is a ruling and $m\geq n+1$.
By following the argument of [W, Theorem 1.1, claim] 
with Theorem 2.13 in [Ko2],
we can prove 

\proclaim{Claim}
Let $X$ be a variety with only log terminal singularities and
$R$ an extremal ray of $X$. 
Let $F$ be an irreducible component of a non-trivial fiber 
of the contraction of $R$.
Assume that for a general point $x\in F$,
there is a rational 
curve $M\subset F$ through $x$ with the following condition:
\roster
\item
its intersection with $-K_X$ is minimal among all rational 
curves in $F$ through $x$.
\item
$X$ has only local complete intersection singularities along $M$
and $M$ is not contained in the singular locus of $X$.
\endroster
Then 
$$\dim F+\dim(\text{locus of}\ R)\geq \dim X+l(R)-1,\tag 1.2.1$$ 
where $l(R)$ is the length of $R$.
Furthermore if the equality holds, the dimension of the deformation
of $M$ which through a fixed point $x$ is $\dim F-1$.
\endproclaim

If $(E_i, -K_X|_{E_i})$ is isomorphic to $(\Bbb P^2, \Cal O_{\Bbb P^2}(2))$,
a general line satisfies the assumption of $M$ in Claim.
So we can use Claim and derive a contradiction to the inequality (1.2.1).
If $(E_i, -K_X|_{E_i})$ is isomorphic to $(\Bbb F_n, C_0+ml)$,
a general ruling satisfies the assumption of $M$ in Claim.
So by using Claim, we obtain the equality in (1.2.1).
But a ruling cannot move if a general point on it is fixed, 
a contradiction to the second part of Claim.
\comment
\item"(3)"
By Theorem 1.19 in [AW2], it is easy to see that
for any two components $E_i$ and $E_j$, $E_i\cap E_j\not= \phi$ and
$E_i\cap E_j$ is a ruling in each component. 
Assume that there are three components $E_i$, $E_j$ and $E_k$ and 
$E_i\cap E_j\not= E_i\cap E_k$. 
By the freeness of $|-K_X|$ ([AW1]), 
we can assume that a general member $C\in |-K_X|$ 
has only Gorenstein terminal singularities.
The restriction of $f$ to $C$ is a flopping contraction.
So $E|_C$ is a tree of rational curves.
But by the assumption, $E_i\cup E_j\cup E_k|_C$ forms a cycle of rational
curves, a contradiction.
Hence all components intersect along a unique curve.
If there are more than $3$ components, $E|_C$ has more than $3$
rational curves which intersects at one point.
But this contradicts to the dual graph
of the minimal resolution of Du Val singularities.
\endcomment
\qed
\enddemo

By this Theorem, the exceptional locus of a Gorenstein terminal
$4$-fold flipping contraction  
is either purely $1$-dimensional or purely $2$-dimensional.
In the former case, we call it a flipping contraction of type $(1,0)$. 
In the latter case, we call it a flipping contraction of type $(2,0)$.

\proclaim{Theorem 1.3}
We consider the object $(*)$.
Let $B$ be a generic hyperplane section through $P$.
Then the strict transform $A:=f^*B$ has only canonical singularities.
\endproclaim
\demo{Proof}
We take a general member $C\in |-K_X|$ and let $D:= f(C)$. 
By the freeness of $|-K_X|$ (Theorem 1.1), 
we can assume that $C$ is Gorenstein terminal.
$D$ is Gorenstein by the Serre-Grothendieck duality 
(cf. [Kaw 2, the Proof of Theorem 8.7]), 
which in turn shows that $D$ is normal and $C\to D$ has then 
only connected fibers 
by the Zariski Main Theorem. 
Hence if $f$ is of type $(1,0)$, $C\to D$ is an isomorphism 
or if $f$ is of type $(2,0)$, $f|_C$ is a flopping contraction.
So in any case $D$ has also Gorenstein terminal singularity at $P$, i.e., 
cDV singularity. 
Then we may assume that $B|_D$ is canonical by replacing $B$ if necessary.
So $A|_C$ must be also normal and canonical
since $A|_C\to B|_D$ is isomorphism if $f$ is of type $(1,0)$ or
$A|_C\to B|_D$ is crepant if $f$ is of type $(2,0)$.
We know that $A$ is canonical along $C|_A$ by the above argument. 
So it suffices to prove that $A$ is canonical outside $C|_A$.
Below argument is inspired by the proof of [Kaw 2, Theorem 8.5].
Let $C'$ be a general member of $|-2K_X|$ and $D':= f(C')$.
Since $K_{B}+(D|_B)|_{D|_B}=K_{D|_B}$ is canonical and $D|_B$
is Cartier on $B$, $K_B + (D|_B)$ is canonical by [Kaw4]
(or [Kaw5, Theorem 1.4]). 
Hence $$K_B + \frac 12 D'|_{B}\  
\text{is also canonical since $D'$ is more general than} \ D. \tag 1.3.1$$
We take the double cover $\tilde{A}\to A$ (resp. $\tilde{B}\to B$)
whose branch locus is $C'|_A$ (resp. $D'|_B$).
Let $g:\tilde{A}\to \tilde{B}$ be the natural morphism.
It is sufficient to prove that $\tilde{A}$ is canonical since $\tilde{A}\to A$
is etale outside $C'|_A$.
By (1.3.1), $\tilde{B}$ is Gorenstein canonical.
So $\tilde{A}$ is also Gorenstein canonical since $g$ is crepant and 
we are done.
\qed
\enddemo
\definition{Remark}
By this Theorem, we see that the object $(*)$ is a very special example of 
a semistable $4$-fold flipping contraction.
(See [C] for the definition of a semistable flipping contraction.) 
\enddefinition

\proclaim{Proposition 1.4}
Consider the situation of Theorem 1.3 and take $A$ and $B$ as there. 
Then 
\roster
\item a general element of $|-K_B|$ has only Du Val singularity at $P$;
\item for any $i$, $E_i$ is not $\Bbb Q$-Cartier divisor in $A$.
\endroster
\endproclaim
\demo{Proof}
\roster
\item $D|_B\in |-K_B|$ in the proof of Theorem 1.3 satisfies $(1)$. 
\item (cf. the argument of [Kac1, 4.3]) 
We assume that for some $i$, $E_i$ is $\Bbb Q$-Cartier.
Assume further that 
$E$ has another component. Let $E_j$ be a component such that
$E_i\cap E_j\not= \phi$. 
Then by the assumption that $E_i$ is $\Bbb Q$-Cartier, $E_i\cap E_j$ is 
$1$-dimensional. 
Then since the Picard numbers of such $E_j$'s and $E_i$ are $1$, 
the union of $E_i$ and $E_j$'s is covered by one extremal ray
in $\overline{NE}(A/B)$. 
For a ruling $m$ of $E_j$ (not contained in $E_i$), $E_i.m>0$.
But for a ruling $l$ in $E_i$, $E_i.l < 0$, a contradiction. 
Hence $E$ is irreducible and $\Bbb Q$-Cartier.
So $B$ has only canonical singularities by [KMM, Lemma 5-1-7].
Note that $B$ is smooth outside $P$ and that  
$|-K_B|$ has a Du Val element through $P$.  
So in fact $B$ is terminal by [St, Section 5]. 
Since $B$ can deform to a $3$-fold with only cDV singularities in $Y$,
$B$ also has only cDV singularity by [Nam, Proposition (3.1)].
In particular $B$ has only hypersurface singularity  
so $Y$ has also only hypersurface singularity, a contradiction.
We establish the proposition.
\endroster
\qed
\enddemo

\proclaim{Theorem 1.5}
Let $V$ be a smooth variety and $X$ a smooth (not necessarily connected)
divisor on $V$.
Let $\pi :V \to S$ be a projective morphism onto a variety $S$
with only connected fibers.
Assume that $K_V +X$ is $\pi$-big for the pair $(V,X)$, i.e.,
$K_V +X$ is $\pi$-big and we can write $l(K_V + X) =A+B$
for a positive integer $l$, a $\pi$-ample divisor $A$ and a $\pi$-effective
divisor $B$ such that 
$\text{Im} (\pi _* \Cal O_V (B) \to \pi _* \Cal O_X (B|_X) )=\not 0$.
Then the natural homomorphism $\pi _* \Cal O_V (m(K_V +X)) \to
\pi _* \Cal O_X (mK_X)$ is surjective for any positive integer $m$.
\endproclaim

\demo{Proof}
See [Si], [Kaw5, 2.2 Theorem A] or [Nak2, Theorem 4.9].
\qed
\enddemo

\comment
\proclaim{Theorem 1.5}
Let $U$ be a smooth $3$-fold and $\mu:U\to V$ a projective bimeromorphic
morphism.  
Let $C$ be the exceptional locus of $\mu$.
Assume that $C$ is $\Bbb P^1$ and $K_U.C > 0$.
Then $\mu$ is target stable, i.e., if there is a small deformation 
$\Cal V\to (S,o)$ of $V$ over $S$, 
then we have a small deformation 
$\Cal U\to \Cal V$ of $\mu: U\to V$ over $(S,o)$.
\endproclaim
\demo{Proof}
See [Ra2, Theorem 3.2].
\qed
\enddemo
\endcomment
 
\proclaim{Proposition 1.6 (H. Laufer)}
Let $S$ be normal Gorenstein surface and 
$f:S\to T$ a projective bimeromorphic morphism to a normal surface $T$.
Let $C$ be the exceptional curve.
Suppose that $C$ is irreducible, isomorphic to $\Bbb P^1$ and $K_S.C = -1$.
Then $f(C)$ is a smooth point of $T$, $S$ has only one singular point
on $C$ which is of type $A_{n-1}$ for some $n\in \Bbb N$. 
Furthermore $C^2= -\frac 1n$.
\endproclaim
\demo{Proof}
See [LS, Theorem 0.1].
\qed
\enddemo

\proclaim{Theorem 1.7(Length of an extremal ray)}
Let $X$ be a variety with only canonical singularities
and $R$ an extremal ray of $X$.
Let $F$ be a $1$-dimensional irreducible component of the fiber
of the contraction of $R$ which contains Gorenstein points of $X$.
Then $K_X.F\geq -1$.
\endproclaim
\demo{Proof}
See [M4, 1.3 and 2.3.2] or [I, Lemma 1].
\qed
\enddemo

\proclaim{Proposition 1.8}
Let $U$ be a $3$-fold with only Gorenstein canonical singularities 
and $(V,P)$ a pair of a $3$-dimensional normal Stein space and a point
in it.
Let $f: U\to V$ be a flipping contraction whose exceptional locus $l$ is
connected.  
Then $l$ is irreducible and
$l \subset \text{Sing}\ U$.
\endproclaim
\demo{Proof}
The irreducibility of $l$ can be proved by the same argument as the first part
of the proof of Theorem 1.3.
Assume that 
$$\text{Sing} \ U  \cap l \ \text{consists of finite points.}\tag 1.8.1$$
Let $g:U'\to U$ be a partial resolution such that $g$ is crepant and
$U'$ has only Gorenstein terminal singularities
(cf. [M3] and [Re2]).
Since $K_{U'}$ is not $f\circ g$-nef, we can find an extremal ray 
$R\in \overline{NE}(U'/V)$.
Let $l'$ be an irreducible curve such that $[l']\in R$. 
Then $l'$ is the strict transform of $l$ by (1.8.1) and
the fact that $K_{U'}$ is $g$-nef. 
So $R$ is a flipping ray. But this contradicts the fact that there is no
flipping contraction from a Gorenstein terminal $3$-fold.
\hfill \qed
\enddemo

\head 2. Proof of the main Theorem  \endhead

\proclaim{Theorem 2.1}
Let $X$ be a $4$-fold with only canonical singularities
and $f\: X\to Y$ be a flipping contraction.
Let $A$ be the pull back of a Cartier divisor $B$ on $Y$.
Assume that $A$ has only canonical singularities.
Then the flip of $f$ exists.
\endproclaim

\demo{Proof}
It suffices to prove that
$\oplus _{m\geq 0} f_* \Cal O (mlK_X)$ is finitely generated,
where $l$ is the minimum positive integer such that $lK_X$ is Cartier.
Let $g:Z\to X$ be a good resolution for the pair $(X,A)$
and $A'$ the strict transform of $A$.
We may take $g$ such that $(\text{excep} g) |_{A'} = \text{excep} (g|_{A'})$.
By Theorem 1.5,
$f_* g_* \Cal O_Z (m(K_Z +A')) \to
f_* g_* \Cal O_{A'} (mK_{A'})$ is surjective for any positive integer $m$.
Note that $(X, A)$ is canonical since $A$ is a Cartier divisor with
only canonical singularities by [Kaw4] (or [Kaw5, Theorem 1.4]).
Hence $f_* g_* \Cal O_Z (ml(K_Z +A')) =f_* \Cal O_X (ml(K_X +A))
=f_* \Cal O_X (mlK_X) \otimes \Cal O_Y (mlB)$ and 
$f_* g_* \Cal O_{A'} (mlK_{A'})=
f_* \Cal O_A (mlK_A)$.
Hence by Nakayama's lemma, it suffices to prove the finite generation of
$\oplus _{m\geq 0} f_* \Cal O (mlK_A)$.
Let $\mu: A_T \to A$ be a small $\Bbb Q$-factorialization,
i.e., $\mu$ is a small projective bimeromorphic morphism such that
$A_T$ is $\Bbb Q$-factorial and has only canonical singularities.  
By running the MMP over $B$ starting from $A_T$ and
taking the canonical model of a minimal model of $A_T$ over $B$,
we obtain ${f^+}':{A^+}' \to B$ such that ${A^+}'$ has only 
canonical singularities and $K_{{A^+}'}$ is ${f^+}'$-ample.
We can easily see that
$f_* \Cal O (mlK_A) ={f^+}' _* \Cal O (mlK_{A^+})$
and $\oplus _{m\geq 0} {f^+}' _* \Cal O (mlK_{A^+})$ is finitely
generated.
Hence we are done.
\qed
\enddemo

\proclaim{Corollary 2.2}
\roster
\item
The flip for a $1$-parameter family of canonical flipping contractions
exists;
\item let $f\: X\to Y$ be a flipping contraction from a Gorenstein
terminal $4$-fold. Then the flip of $f$ exists.
\endroster
\endproclaim

\demo{Proof}
(1) is clear.
By Theorem 1.3, we can apply Theorem 2.1 for (2).
\qed
\enddemo

From now on we consider the object as in $(*)$.
We will use the notation as in the proof of Theorem 1.3 freely.
We will denote by $f^+\: X^+\to Y$ the flip of $f$ and 
by $E^+$ the exceptional locus of $f^+$. 
Furthermore we will denote with $+$ the strict transform of a divisor of $X$ 
on $X^+$. 

\proclaim{Corollary 2.3} 
\roster
\item $\dim E^+ =1$ and $\dim E =2$;
\item $A^+$ and $X^+$ have only Gorenstein terminal singularities; 
\item $f^+|_{C^+}$ is the $-K_X|_C$-flop for $f|_C$.
In particular
$\# \{\text{components of} \ E \}=
\# \{\text{components of} \ E|_C \}=
\# \{\text{components of} \ E^+ \}$;
\item assume that $(X, E)$ is $\Bbb Q$-factorial.
Then $E$ and $E^+$ are irreducible. 
\item assume that $E$ is irreducible
(and hence isomorphic to $\Bbb F_{n,0}$ for some natural
number $n$).
Then
there exists a Weil divisor $H$ on $X$ such that $-K_X\sim nH$. 
By using $C$ (as in the proof of Theorem 1.3), 
define ring structures to
$$\bigoplus^{n-1}_{j = 0} \Cal O_X(-jH) \
\text{and} \ \bigoplus^{n-1}_{j = 0} \Cal O_Y(-jf(H))$$
and set $$\tilde{X}:= \pmb {Specan} \bigoplus^{n-1}_{j = 0} \Cal O_X(-jH) 
\ \text{and}  
\ \tilde{Y}:= \pmb {Specan} \bigoplus^{n-1}_{j = 0} \Cal O_X(-jf(H)).$$
Let $\tilde{E}\subset \tilde{X}$ be the pull back of $E$
by the natural morphism $\tilde{X}\to X$.
Then the natural morphism 
$\tilde{f}:\tilde{X}\to \tilde{Y}$ is a flipping contraction
which satisfies the same assumption as $f$, and 
$\tilde{E}$ is the exceptional locus of $\tilde{f}$ and 
is isomorphic to $\Bbb P^2$. 
\endroster
\endproclaim

\demo{Proof}
\roster
\item 
We saw in the proof of Theorem 1.3 
that $P$ is a Gorenstein terminal singularity of $D$. 
On the other hand $E^+ \subset C^+$ and $f^+ |_{C^+}$ is crepant.
Hence $\dim E^+ =1$.
Furthermore by [KMM, Lemma 5-1-17], $\dim E =2$;
\item 
Since $P$ is a canonical singularity of $D\cap B$
and $f|_{C^+ \cap A^+}$ is crepant,
$C^+ \cap A^+$ has only canonical singularity.
Hence by $\text{Sing} A^+ \subset E^+$ and [St, Section 5],
$A^+$ has only terminal singularities. 
Furthermore by the argument as in the proof of Proposition 1.4 (2),
we see that $A^+$ is Gorenstein and so is $X^+$;
\item
since $A^+ \cap C^+$ has only canonical singularity
and $\text{Sing} C^+ \subset E^+$,
$C^+$ has only terminal singularities.
Note that $-K_X |_{C}$ is $f|_C$-ample and 
$-K_{X^+} |_{C^+}$ is ${f^+}|_{C^+}$-negative.
Hence by the uniqueness of the $-K_X|_C$-flop, 
$f^+|_{C^+}$ is the $-K_X|_C$-flop for $f|_C$.
\item 
By $\Bbb Q$-factoriality of $(X,E)$,
$\rho (X^+ /Y, E^+)=1$ whence $E^+$ must be irreducible.
Hence $E$ is also irreducible by (3).
\item
\comment
Since $C$ is smooth and
the analytic structure along the exceptional curves is unchanged
by the $-K_X|_C$-flop (cf.[Ko1, Theorem 2.4]),
$C^+$ is also smooth.
Let $E'$ be the exceptional curve of $f|_C$
and $C'\in |-K_X|_C|$.
Then we have $C'.E' = n$ by Theorem 1.2. 
Let ${C'}^+\in |-K_{X^+}|_{C^+}|$ be the strict transform of $C'$.
\endcomment
Since $C\dashrightarrow C^+$ is a terminal flop, we have
${C'}^+.E^+ = - n$.
Let $H^+$ be a Cartier divisor on $X^+$ such that $H^+.E^+ = -1$
and $H \subset X$ be the strict transform of $H^+$.
Then $K_{X^+} + nH^+$ is linearly $f^+$-trivial 
since we consider locally analytically along $E^+$.
Hence $K_X + nH$ is linearly $f$-trivial since the linear triviality is 
preserved by an anti-flip.
Let $\t{X}$ and $\t{Y}$ are as in the statement of (5). 
We check that the natural morphism $\tilde{f}:\tilde{X} \to \tilde{Y}$
and the pull backs $\tilde{A}$ of $A$, $\tilde{B}$ of $B$
and $\tilde{E}$ of $E$ 
satisfy the same assumption as $f$, $A$, $B$ and $E$ 
and we prove that $\tilde{E}$ is $\Bbb P^2$. 
Let $\pi:\tilde{X}\to X$ be the covering morphism and $\tilde{C}:=
(\pi^*(C))_{\text{red}}$.
Note that $n\tilde{C}=\pi^*(C)$, $\tilde{C}\simeq C$ and 
$\tilde{C}$ is a Cartier divisor since $C$ is contained in the branch locus.
Then by the ramification formula $K_{\tilde{X}}=\pi^*K_{X}+(n-1)\tilde{C}$,
$\tilde{C}\in |-K_{\tilde{X}}|$. Since $\tilde{C}$ is a Cartier divisor,
we see that $\tilde{X}$ is Gorenstein. 
We also know that $\tilde{X}$ is terminal
since codimension $1$ ramification locus $\tilde{C}$ of $\pi$ has only 
terminal singularities.
The rest are clear except that $\tilde{E}\simeq \Bbb P^2$. 
The restriction of $\pi$ to $\tilde{E}$ is
$\pi|_{\tilde{E}}:
\tilde{E}=\pmb {Specan} \bigoplus^{n-1}_{j = 0} \Cal O_E(-jl)\to E$, 
where $l$ is a ruling of $E$. (Note that $H|_E\sim l$.)
So it coincides with the quotient $\Bbb P^2\to \Bbb F_{n,0}$ by the action
of $\Bbb Z_n$, $(X:Y:Z)\to (\eta X:\eta Y:Z)$, 
where $X,Y\ \text{and} \ Z$ is the homogeneous coordinate of $\Bbb P^2$
and $\eta$ is a primitive $n$-th root of unity. So $\tilde{E}$ is $\Bbb P^2$. 
\endroster
\qed
\enddemo

We can classify $f$ and $f^+$ with additional assumptions as follows:

\proclaim{Corollary 2.4}
Assume that $A$ (as in Theorem 1.3) 
has only isolated singularities
and $E$ is irreducible (and hence isomorphic to $\Bbb F_{n,0}$ for some natural
number $n$).
Then
\roster
\item $A$ is singular only at the vertex $v$ of $E$ 
(if $n=1$, the vertex means a point on $E$).
Near $v$, $(v\in E\subset A\subset X)$ 
is analytically isomorphic to $(o\in (x=z=t=0)\subset (xy+zw =t= 0) \subset
(xy+zw+t^k=0)) \ \text{in} \ \Bbb C ^5 /{\Bbb Z_n (1,-1,1,-1,0)}$;
\item $X^+$ is smooth, $E^+$ is $\Bbb P^1$ and $\Cal N_{E^+/X^+}\simeq
\Cal O(-1)\oplus \Cal O(-1)\oplus \Cal O(-n)
\ \text{or} \ 
\Cal O\oplus \Cal O(-2)\oplus \Cal O(-n).$
Furthermore the former case occurs if and only if $X$ has only 
$\frac 1n (1,-1,1,-1)$ singularity at $v$;
\endroster
\endproclaim

\demo{Proof}
We will prove (1).
Consider the covering as in Corollary 2.3 (5).
We will use the notation as in its proof.
Then $\t{A}$ has only isolated singularities since $C\cap A$ is smooth.
Let $q:\t{A_q}\to \t{A}$ be a small morphism such that the inverse image 
$\t{E_q}$
of $\t{E}$ is $q$-anti-ample (i.e., 
$\t{A_q}:= \pmb{Projan}\bigoplus^{\infty}_{m=0} \Cal O_{\t{A}}(-m\t{E})$ 
and $q$ is the natural projection.)
We can take such a small morphism by [Kaw2, Theorem 6.1].
Since $\t{E}$ is not $\Bbb Q$-Cartier by proposition 1.4, 
$\t{A_q}$ is not isomorphic to $\t{A}$. 
Let $\Phi:\t{A_q}\to \t{A^+}$ be the contraction of an extremal ray in 
$\overline{NE}(\t{A_q}/\t{B})$ and $g^+:\t{A^+}\to \t{B}$ the natural morphism.
We obtain the following diagram:
$$\matrix 
\ & {\t{A_q}} & \  \\
{q\swarrow} & \ & {\searrow\Phi} \\
\t{A}  &\ & {\t{A^+}} \\
{\t{f}|_{\t{A}}}\searrow & \ & {\swarrow g^+}\\
\ & \t{B} & \ & . \\
\endmatrix $$

\proclaim{Claim 2.5}
$\Phi$ is a divisorial contraction which contracts $\t{E_q}$ to a curve.
\endproclaim
\demo{Proof}
Since $q$ is not an isomorphism and $-\t{E_q}$ is $q$-ample,
$\t{E_q}$ contains all $q$-exceptional curves.
If $\Phi$ is a divisorial contraction which contracts $\t{E_q}$ to a point,
such $q$-exceptional curves are contracted by $\Phi$.
But this is absurd since $K_{\t{A_q}}$ is $q$-trivial but $\Phi$-negative. 
If $\Phi$ is a flipping contraction, then the flipping curve $m$
is contained in the curve singularity of $\t{A_q}$ by Proposition 1.8.
So by the assumption that $\t{A}$ has only isolated singularities,
$m$ must be contained in the $q$-exceptional curve, a contradiction.
\qed
\enddemo

Let $\t{E^+}$ be the curve $\Phi(\t{E_q})$.

\proclaim{Claim 2.6}
\roster
\item 
$\t{A^+}$ is smooth along $\t{E^+}$. 
$\t{E^+}\simeq \Bbb P^1$ and $\t{E_q}\simeq \Bbb F_1$;
\item 
$g^+:\t{A^+}\to \t{B}$ is isomorphic to the restriction of $\t{f^+}$
to the strict transform of $\t{A}$ on $\t{X^+}$.
(Hence we will denote by $\t{A^+}$ 
the strict transform of $\t{A}$ on $\t{X^+}$.) 
$\t{A_q}$ is smooth and $\tilde{E_q}.M=-1$, 
where $M$ of $\t{E_q}$ is the negative section.
\endroster
\endproclaim

\demo{Proof}
\roster
\item
By Theorem 1.1, 
$|-K_{\t{A_q}}|$ is free near the fiber over any point $Q$ of $\t{E^+}$,
so we can take a smooth member $D\in |-K_{\t{A_q}}|$ 
near the fiber since $\t{A}$ has only canonical singularities. 
Since $D$ maps isomorphically to $\Phi(D)\in |-K_{\t{A^+}}|$ 
(cf.[Kaw2, the Proof of Theorem 8.7]),
we see that there is a smooth member of $|-K_{\t{A^+}}|$ through $Q$.
Note that $Q$ is a canonical singularity of $\t{A^+}$.
By these, we can see that $Q$ is a smooth point of $\t{A^+}$ as follows:

It is sufficient to prove that $K_{\t{A^+}}$ is Cartier at $Q$.
Assume the contrary. Let $\pi: \o{A}^+ \to \t{A^+}$ be the index $1$ cover
for $K_{\t{A^+}}$ near $Q$. 
Then $\o{A}^+$ is Gorenstein canonical at $\pi^{-1}(Q)$.
Since $\pi$ is ramified only at $Q$ and $\Phi(D)$ is smooth,
$$\pi^{-1}\Phi(D) \ \text{has at least} \  2 
\ \text{components and they intersect mutually only at} \
\pi^{-1}(Q).\tag 2.6.1$$
Furthermore they are all smooth. In particular
$\pi^{-1}\Phi(D)$ satisfies $R_1$ condition.
On the other hand $\pi^{-1}\Phi(D)$ satisfies $S_2$ condition
since this is a Cartier divisor of a canonical singularity.
Hence $\pi^{-1}\Phi(D)$ is normal by the Serre's criterion.
But this is a contradiction to (2.6.1).

Since $\t{B}$ has only rational singularities 
and $\t{E^+}$ is an irreducible curve,
$\t{E^+}$ must be $\Bbb P^1$.
Since a general fiber $n$ of $\Phi$ is irreducible and reduced and 
$-K_{\t{A_q}}.n = 1$ (Theorem 1.7), any fiber is irreducible and reduced.
So $\t{E_q}$ is $\Bbb F_1$.

\item
Let $Q$ be any point on $\t{E^+}$, $G$ a general (smooth) 
hyperplane section of $\t{A^+}$
through $Q$ such that $\t{A^+}|_{\t{E^+}}$ is one point 
and $F$ the pull back of $G$ (we consider analytically locally
near $Q$).
Then $F$ is normal.
In fact, since the fiber $(\t{E_q})_Q$ over $Q$ of $\Phi$ is not contained 
in the singular locus of $\t{A_q}$, 
$\t{E_q}$ is generically Cartier in $\t{A_q}$ along $(\t{E_q})_Q$, 
which in turn shows $(\t{E_q})_Q$ is generically Cartier divisor on $F$. 
Since $(\t{E_q})_Q$ is smooth, $F$ is generically smooth along 
$(\t{E_q})_Q$, 
i.e., $F$ is normal. 
Furthermore we have $K_F.(\t{E_q})_Q = -1$ and $F$ is 
Gorenstein.
So we know by Proposition 1.6 that $F$ has only one $A_{m-1}$ singularity
for some integer $m$ and $((\t{E_q})_Q)^2_F = -\frac 1m$.
On the other hand, $((\t{E_q})_Q)^2_F = (\t{E_q}.(\t{E_q})_Q)_{\t{A_q}}$ 
and the value of
the right side of this equality is independent of $Q$,  
so $m$ is also independent of $Q$. 
Hence we find that $\t{A_q}$ has the locally trivial 
$cA_{m-1}$ curve singularity
along $M$ and outside $M$, $\t{A_q}$ is smooth. 
By $((\t{E_q})_Q)^2_F = -\frac 1m$,
we obtain the subadjunction formula
$K_{\t{E_q}} + \frac {m-1}{m}M = K_{\t{A_q}} + \t{E_q}|_{\t{E_q}}$
and $K_{\t{A_q}} = \Phi^*K_{\t{A^+}} + m\t{E_q}$.
Intersecting these with $M$, we can see that 
$K_{\t{A^+}} .\t{E^+} =2m-1$.
Remark that $\t{E_q}.M$ is negative since $\t{E_q}$ is $q$-anti-ample.
So $K_{\t{A^+}}.\t{E^+}$ is positive and
hence $g^+ :\t{A^+}\to \t{B}$ is the canonical model of $\t{f}|_{\t{A}}$.
On the other hand 
the restriction of $\t{f^+}$
to the strict transform of $\t{A}$ on $\t{X^+}$ is also 
the canonical model of $\t{f}|_{\t{A}}$ by Corollary 2.3.
Hence by the uniqueness of the canonical model, they are isomorphic.
So we have $2m-1=-1$, i.e., $\t{A_q}$ is smooth also along $M$.
\endroster
\qed
\enddemo

By considering the normal bundle sequence 
$$0\to \Cal N_{M/{\t{E_q}}} \to  
\Cal N_{M/{\t{A_q}}} 
\to \Cal N_{\t{E_q}/\t{A_q}}|_M \to 0,$$
we see that
$\Cal N_{M/\t{A_q}}\simeq 
\Cal O(-1)\oplus \Cal O(-1)$. 
So $M$ is contracted to an ordinary
double point by $q$. Denote this point by $\t{v}(\in \t{A})$.
We note that $\t{X}$ is singular at worst only at $\t{v}$
since so is $\t{A}$ and $\t{X}^+$ is smooth.
Then $\t{v}$ is the unique isolated ramification point of $\pi$ 
and hence $\t{v} = \pi^{-1}(v)$ (Recall that $v$ is the vertex of $E$.)
So we can write locally analytically 
$(\t{v}\in \t{E}\subset \t{A}\subset \t{X}) 
\simeq (o\in (x = z = t = 0)\subset
(xy + zw = t = 0)\subset (xy+zw+t^k=0))$, where
$x, y, z, w$ are the semi-invariant coordinates	
and $xy + zw$ is semi-invariant with respect to the action of $\Bbb Z_n$.
When we restrict the action to $\t{E}$, the action is $(y, w)\to 
(\eta y, \eta w)$, where $\eta$ is a primitive $n$-th root of unity by the
explicit description of $\pi|_{\t{E}}$.
Hence the action is $(x, y, z, w)\to (\eta^a x, \eta y, \eta^a z, \eta w)$, 
where $a$ is an integer.
By the necessary condition for the quotient to be canonical ([M3, Theorem 2]), 
$a$ must be $-1$. This is also sufficient.

Next we will prove (2).
To determine the normal bundle $\Cal N_{E^+/X^+}$,
we consider the normal bundle sequence
$$0\to \Cal N_{E^+/C^+} \to \Cal N_{E^+/X^+} 
\to\Cal N_{C^+/X^+}|_{E^+} \to 0.\tag 2.4.1$$
Since $A|_C$ is smooth, we see that 
$\Cal N_{E|_C/C}= \Cal O\oplus \Cal O(-2)\ \text{or} \ 
\Cal O(-1)\oplus \Cal O(-1)$ by using the normal bundle sequence
$$0\to \Cal N_{{E|_C}/{A|_C}} \to \Cal N_{{E|_C}/C} 
\to\Cal N_{A|_C/C}|_{E|_C} \to 0.$$
Since $C\dashrightarrow C^+$ is the flop, we have 
$\Cal N_{E^+/C^+}\simeq \Cal N_{{E|_C}/C}$.
On the other hand, $\Cal N_{C^+/X^+}|_{E^+} = \Cal O(-n)$, so
the sequence (2.4.1) is split. Hence we obtained the first part of $(2)$.

To prove the second part of (2), we consider the covering described in 
Corollary 2.3 (5).
We use the notation there. Recall that $\tilde{C}\simeq C$ and $\tilde{C}\in
|-K_{\tilde{X}}|$. By the argument above together with this,
we see that
$$\Cal N_{E^+/X^+}\simeq
\Cal O(-1)\oplus \Cal O(-1)\oplus \Cal O(-n) \ \text{if and only if} \ 
\Cal N_{{\tilde{E}|_{\tilde{C}}}/{\tilde{C}}}\simeq
\Cal O(-1)\oplus \Cal O(-1). \tag 2.4.2$$
So 'if' part of (2) follows from Kawamata's determination of a flipping
contraction from a smooth $4$-fold (see Theorem 0.1) and (1).
Finally we prove 'only if' part. Assume that
$\Cal N_{E^+/X^+}\simeq
\Cal O(-1)\oplus \Cal O(-1)\oplus \Cal O(-n)$. 
Then by (2.4.2), 
$\Cal N_{{\tilde{E}|_{\tilde{C}}}/{\tilde{C}}}\simeq
\Cal O(-1)\oplus \Cal O(-1)$. 
Then locally analytically there is a smooth surface $\tilde{S}$ such that
$\tilde{S}\subset \tilde{C}$ and $\tilde{S}.(\tilde{E}|_{\tilde{C}})=-1$ 
(note that $\tilde{S}\in |K_{\tilde{X}}|_{\tilde{C}}|$.)
Let $\tilde{S}^+\in |K_{\tilde{X}^+}|_{\tilde{C}^+}|$ be 
the strict transform of $\tilde{S}$ on $\tilde{C}^+$.
Consider the exact sequence 
$$0\to \Cal O_{\tilde{X}^+}(2K_{\tilde{X}^+})\to 
\Cal O_{\tilde{X}^+}(K_{\tilde{X}^+}) \to
\Cal O_{\tilde{C}^+}(K_{\tilde{X}^+}) \to 0.$$ 
By the Kodaira-Kawamata-Viewheg vanishing theorem, 
$H^1(\tilde{X}^+, \Cal O_{\tilde{X}^+}(2K_{\tilde{X}^+}))=0$. 
So there is an element $\tilde{V}^+\in |K_{\tilde{X}^+}|$ such that 
$\tilde{V}^+|_{\tilde{C}^+}= \tilde{S}^+$.
Let $\tilde{V}\in |K_{\tilde{X}}|$ be the strict transform of $\tilde{V}^+$.
Then $\tilde{V}|_{\tilde{C}}=\tilde{S}$. We claim that $\tilde{V}$ is smooth. 
This implies $\tilde{X}$ is also smooth,
which completes the proof of the 'only if' part whence (1).
First we note that $\tilde{V}$ is normal 
since $\tilde{V}|_{\tilde{C}}$ is smooth.
Let $x$ be any point of $\tilde{V}$ and 
$\tilde{C}_x$ a normal general member of
$|-K_{\tilde {X}}|_{\tilde{V}}|$.
Let $\tilde{E}_x:= \tilde{E}|_{\tilde{C}_x}$. 
$\tilde{E}_x$ is the exceptional curve of $\tilde{f}|_{\tilde{C}_x}$ and
$K_{\tilde{C}_x}.\tilde{E}_x=-1$.
Hence by Proposition 1.6, $\tilde{C}_x$ has only one singular point on 
$\tilde{E}_x$
which is of type $A_{m-1}$ for some $m$ and 
$(\tilde{E}_x)^2_{\tilde{C}_x}=-\frac1m$.  
Since $(\tilde{E}_x)^2_{\tilde{C}_x}= 
(\tilde{E}^2.{\tilde{C}_x})_{\tilde{V}}$, $m$ is independent of $x$.  
Hence $(\tilde{E}_x)^2_{\tilde{C}_x}= 
(\tilde{E}|_{\tilde{S}})^2 = -1$ and $m$ is $1$, i.e., 
$\tilde{C}_x$ is smooth. 
Consequently $\tilde{V}$ is found to be smooth at any point $x$ and 
we are done.
Now we finished the proof of Corollary 2.4.
\qed
\enddemo

\comment
\definition{Remark 2.4}
Let $f:X\to (Y,P)$ be a flipping contraction from a Gorenstein terminal
$4$-fold. We use the notation of Theorem 1.3 and Proposition 1.4.
By Proposition 1.5, a general member of $|-K_B|$ has only Du Val
singularity at $P$. So $\bigoplus^{\infty}_{j = 0} \Cal O_X(-jK_B)$ and
$\bigoplus^{\infty}_{j = 0} \Cal O_X(jK_B)$ are finitely generated.
Set $A^+:= \pmb{Projan} \bigoplus^{\infty}_{j = 0} \Cal O_X(jK_B)$
and let $g^+: A^+\to B$ be the natural morphism.
Then $A^+$ has only terminal singularities,
$K_{B^+}$ is $g^+$-ample and $g^+$ is a small morphism 
(see [KM, Theorem 3.1]).
So if we can generalize Theorem 1.5, we may construct the flip of $f$
by deforming $g^+$.
\enddefinition
\endcomment

\head 3. Some examples \endhead

We construct examples of flipping contractions 
from Gorenstein terminal $4$-folds. 

\definition{Example 3.1 (Toric example)}
Let $\pmb e_i$ be the vector 
$(0,.., \overset{\overset{i}\to {\lor}} \to {1},.., 0)$ in $\Bbb R^4$
for $i = 1, 2, 3$, $\pmb e_4 = (-1, -1, n-1, n)$
and $\pmb e_5 = (0, 0, -1, -1)$.
Let $C_i$ be the cone $<\pmb e_1, \pmb e_2, .., \check{\pmb e_i},.., \pmb e_5>$
for $i\geq 0$ and $C_0$ 
the cone $<\pmb e_1, \pmb e_2,.., \pmb e_5>$.
We denote the toric variety associated to the fan $*$ by $V(*)$. 
Set $X:= V(C_3\cup C_4\cup C_5)$, $X^+ := V(C_1\cup C_2)$
and $Y:= V(C_0)$.
Let $f: X\to Y$ and $f^+: X^+\to Y$ be the natural morphisms.
Then it is easy to check that they define a flipping diagram.
(See [Re].)
\enddefinition

\definition{Example 3.2 (Y. Kachi, M. Gross)}
For the above example we can easily find $A$ (as in the main theorem)
with only isolated canonical singularity as determined in Corollary 2.4.
We can consider that $X$ is locally 
a $1$-parameter family of $A$ over the unit disk $\Delta(t)$.   
Take the cyclic coverings $\hat{X}\to X$, $\hat{Y}\to Y$ and
$\hat{X^+}\to X^+$ associated to the cyclic covering $\Delta(s)\to \Delta(t)$
defined by $t = s^m$.
Then the natural morphisms $\hat{X}\to \hat{Y}$ and $\hat{X^+}\to \hat{Y}$
give a flipping diagram.
\enddefinition

\definition{Example 3.3} 
For the example 3.1 with $n=1$, we can find $A$ 
whose singularity is the curve singularity of generically $cA_1$ type 
along a line of $\Bbb P^2$.
For this $A$, we make the similar construction to Example 3.2. 
We obtain a flipping contraction from a Gorenstein terminal $4$-fold
which has a $1$-dimensional singular locus.
Furthermore if we take $q:A_q\to A$ as in the proof of the main theorem
for this $A$, the first extremal contraction of $A_q$ over $B$ is
a flipping contraction and after the flip, we can contract the strict
transform of $E$ to a Gorenstein terminal point.
\enddefinition

\enddocument